\documentclass[11pt]{article}
\usepackage{times}
\usepackage{latexsym}
\usepackage{epsfig,enumerate,amsmath,amsfonts,amssymb,setspace}
\usepackage{indentfirst}
\usepackage{wrapfig,lipsum}
\usepackage{epsfig}
\usepackage{graphicx}
\usepackage{amsthm}
\usepackage{algorithm}
\usepackage{algorithmic}
\usepackage[T1]{fontenc}
\usepackage{authblk}
\usepackage[algo2e,ruled,vlined]{algorithm2e}
\newtheorem{theorem}{Theorem}[section]
\newtheorem{lemma}{Lemma}[section]

\newtheorem{obs}{Observation}[section]
\newtheorem{claim}{Claim}[section]

\newcommand{\pf}{{\bf Proof : }}

\title{A linear-time algorithm for semitotal domination in strongly chordal graphs}
\author[1]{Vikash Tripathi\thanks{2017maz0005@iitrpr.ac.in}}
\author[1]{Arti Pandey\thanks{arti@iitrpr.ac.in}}
\author[2]{Anil Maheshwari \thanks{anil@scs.carleton.ca}}
\affil[1]{Department of Mathematics, Indian Institute of Technology Ropar, Punjab, India.}
\affil[2]{School of Computer Science, Carleton University, Ottawa, Canada.}

\date{}
\begin{document}
\maketitle
\begin{abstract}
In a graph $G=(V,E)$ with no isolated vertex, a  dominating set $D \subseteq V$, is called a \emph{semitotal dominating} set if for every vertex $u \in D$ there is another vertex $v \in D$, such that distance between $u$ and $v$ is at most two in $G$. Given a graph $G=(V,E)$ without isolated vertices, the \textsc{Minimum Semitotal Domination} problem is to find a minimum cardinality semitotal dominating set of $G$. The \emph{semitotal domination number}, denoted by $\gamma_{t2}(G)$, is the minimum cardinality of a semitotal dominating set of $G$. The decision version of the problem remains NP-complete even when restricted to chordal graphs, chordal bipartite graphs, and planar graphs. Galby et al. in \cite{semi1} proved that the problem can be solved in polynomial time for bounded MIM-width graphs which includes many well known graph classes, but left the complexity of the problem in strongly chordal graphs unresolved. Henning and Pandey in \cite{semi2} also asked to resolve the complexity status of the problem in strongly chordal graphs. In this paper, we resolve the complexity of the problem in strongly chordal graphs by designing a linear-time algorithm for the problem.
\end{abstract} 

\section{Introduction}
\label{intro}

A \emph{dominating set} in a graph $G=(V,E)$, is a set $D \subseteq V$, such that any vertex not in $D$ is adjacent to a vertex in $D$. The minimum size of a dominating set is called \emph{domination number}, denoted by $\gamma(G)$. The \textsc{Minimum Domination} problem involves computing a minimum cardinality dominating set of a graph $G$. The domination number is one of the most studied parameter in the graph theory. A thorough treatment and detailed study on domination can be found in the books \cite{dbook4,dbook3,dbook1,dbook2}. Due to numerous applications in the real world problems, many researchers introduced several variations of domination by imposing one or more additional conditions on dominating set. One of the most important variation of domination is total domination.

In a graph $G=(V,E)$, without isolated vertices, a dominating set $D \subseteq V$ is called a \emph{total dominating set} (TD-set in short), if $G[D]$, the graph induced by $D$ in $G$ has no isolated vertex. The \emph{total domination number}, denoted by $\gamma_{t}(G)$, is the cardinality of a minimum total dominating set of $G$. The \textsc{Minimum Total Domination} problem requires to compute a total dominating set of a graph $G$ with no isolated vertex, of size $\gamma_{t}(G)$. See \cite{td1,tdbook} for the detailed results on total domination.
%The details of result on total domination can be found in the book \cite{tdbook} and in survey paper \cite{td1}.

Goddard, Henning, and McPillan, introduced a relaxed notion of total domination, called \emph{semitotal domination} in \cite{semi3} and further studied in \cite{semi1,semi8,semi10,semi7,semi5,semi12,semi4,semi14,semi2,semi11,semi13,semi9,semi6}, from both algorithmic and combinatorial point of view. In a graph $G$ with no isolated vertices, a \emph{semitotal dominating set}(in short, semi-TD-set) is a dominating set $D \subseteq V$ such that for every vertex $u \in D$, there is another vertex $v \in D$, such that the distance between $u$ and $v$ is at most two in $G$. The \emph{semitotal domination number}, denoted by $\gamma_{t2}(G)$, is the cardinality of a minimum semi-TD-set of $G$. It follows directly from definitions that every total dominating set is a semitotal dominating set. Hence, for a graph $G$ with no isolated vertices, we have the following relation between the three parameters:
\vspace{-0.1cm}
\[
\gamma(G) \leq \gamma_{t2}(G) \leq \gamma_{t}(G).
\]

Therefore, the semitotal domination number is squeezed between two important parameters, domination number and total domination number. The minimum semitotal domination problem and and its decision version are defined as follows:
\smallskip

\noindent\underline{\textsc{Minimum Semitotal Domination} problem}
\\
[-12pt]
\begin{enumerate}
  \item[] \textbf{Instance}: A graph $G=(V,E)$ with no isolated vertices.
  \item[] \textbf{Solution}: A Semi-TD-set $D$ of $G$.
  \item[] \textbf{Measure}: Cardinality of the set $D$.
\end{enumerate}

\noindent\underline{\textsc{Semitotal Domination Decision} problem}
\\
[-13pt]
\begin{enumerate}
  \item[] \textbf{Instance}: A graph $G=(V,E)$ and a positive integer $k \le |V|$.
  \item[] \textbf{Question}: Does there exist a Semi-TD-set $D$ in $G$ such that $|D| \le k$?
\end{enumerate}

\begin{figure}[ht]
\begin{center}
\includegraphics[width=0.95\textwidth]{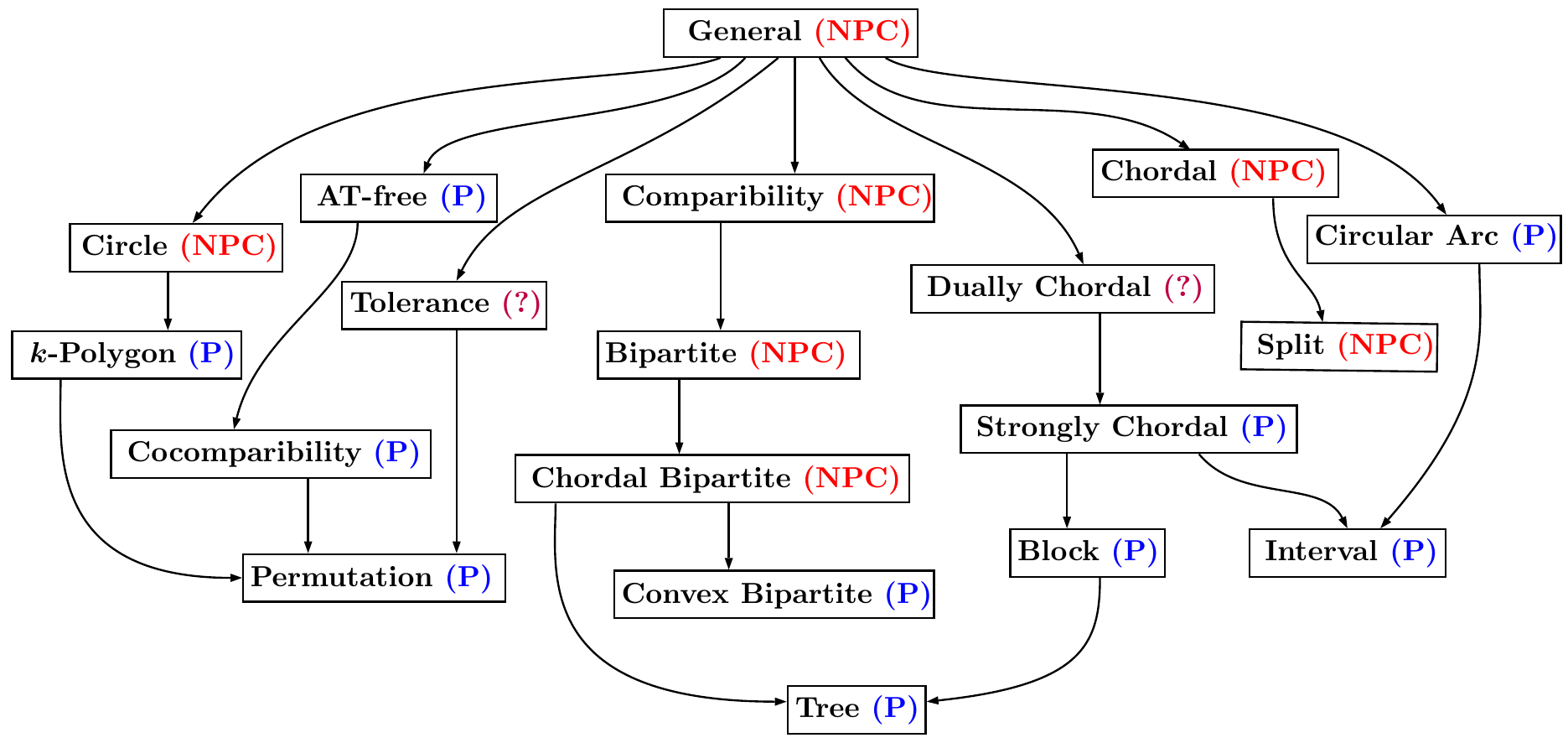}
\caption{Complexity status of \textsc{Minimum Semitotal Domination} problem in some well known graph classes.}
\label{fig:1}
\end{center}
\end{figure}

The \textsc{Semitotal Domination Decision} problem in NP-complete~\cite{semi3} for general graphs. The problem remains NP-complete, even when restricted to chordal graphs, chordal bipartite graphs, and planar graphs~\cite{semi2}. On positive side, we have polynomial-time algorithms to compute a minimum cardinality semi-TD-set in trees \cite{semi3}, interval graphs \cite{semi2,semi13} and block graphs \cite{semi14}. Galby et al. \cite{semi1} proved that, a minimum semi-TD-set can be computed in polynomial-time in bounded MIM-width graphs, which includes many important graph classes.
The complexity status of the problem in some well known graph classes is shown in Fig.~\ref{fig:1}. In the figure, P stands for polynomial-time and NPC stands for NP-complete. The complexity status of the problem in graph classes with question mark is still unknown.

Henning and Pandey studied the approximation hardness of \textsc{Minimum Semitotal Domination} problem \cite{semi2}. They proved that the problem can not be approximated within $(1-\epsilon)$ln$(|V|)$ for any $\epsilon >0$, unless $\text{NP} \subseteq \text{DTIME}(|V|^{O(\text{loglog}(|V|)})$. On other side, they proved that the \textsc{Minimum Semitotal Domination} problem is in the class log-APX. They also proved that the problem is APX-complete for the bipartite graphs with maximum degree $4$.

Although, we have noticed that the semitotal domination number is squeezed between domination number and total domination number. But, \textsc{Minimum Semitotal Domination} problem and \textsc{Minimum Total Domination} problem differs in complexity, see~\cite{semi2}. Indeed, the \textsc{Minimum Total Domination} problem is polynomial-time solvable in chordal bipartite graphs but decision version of \textsc{Minimum Semitotal Domination} is NP-complete for chordal bipartite graphs. Further, Galby et. al. proved that it is NP-hard to decide $\gamma_{t2}(2) =\gamma_{t}(G)$, even when $G$ is a planar graph with degree at most $4$, see~\cite{semi1}.

As the MIM-width of strongly chordal graphs is unbounded, the complexity of the problem was left open in strongly chordal graphs by Galby et al. \cite{semi1}. Henning and Pandey in~\cite{semi2}, also asked to find the complexity status of the problem in strongly chordal graphs. In this paper, we prove that the \textsc{Minimum Semitotal Domination} problem can be solved in linear-time in strongly chordal graphs.

The further structure of the paper is as follows. In Section~\ref{not}, we discuss some notations and definitions. In Section~\ref{sec:1}, we discuss strongly chordal graphs and their properties. In Section~\ref{sec:2}, we design a linear-time algorithm to compute a semi-TD-set in strongly chordal graphs. Finally, Section~\ref{con}, concludes the paper.

\section{Preliminaries}
\label{not}
Let $G=(V,E)$ be a simple graph, where $V=V(G)$ and $E=E(G)$. Two distinct vertices $u,v \in E(G)$, said to be \emph{adjacent} if $uv \in E(G)$. For a vertex $v \in V$, the set $N_{G}(v)= \{u \in V \mid uv \in E(G)\}$ denotes the \emph{open neighbourhood} of $v$ in $G$ and the set $N_{G}[v]=N_{G}(v) \cup \{v\}$ denotes the \emph{closed neighbourhood} of $v$ in $G$. A \emph{path} $P=v_0v_1\ldots v_k$, is a sequence of distinct vertices, such that $v_{i-1}v_{i} \in E(G)$, where $1 \leq i \leq k$ and $k \geq2$. Such a path, is called, a path between $v_0$ and $v_k$. We denote $V(P) =\{v_0,v_1, \ldots v_k\}$.  The length of the path $P$ is $|V(P)|-1$. The distance between two distinct vertices $u,v \in V(G)$, denoted by $d_{G}(u,v)$, is the length of the shortest path between $u$ and $v$ in $G$. Further, we call $u$, a \emph{distance two neighbour} of $v$, if $d_{G}(u,v) \leq 2$.

A path $P=v_0v_1\ldots v_k$ with an additional condition that, $v_0v_k \in E(G)$ is known as a \emph{cycle} on $k$-vertices, denoted by $C_k$. In a cycle $C_k$, where $k\geq 4$, a \emph{chord} is an edge joining two non-consecutive vertices of $C_k$. A graph $G$ is called \emph{chordal}, if any cycle of length at least $4$ in $G$, has a chord.

Let $|V(G)| = n$, and $\beta=(v_1,v_2, \ldots,v_n)$ be any ordering of the vertex set $V(G)$. For a vertex $v_i$ in the ordering $\beta$, we define the sets $N_{i}(v_i) = \{v_j \mid j >i ~\text{and} ~v_iv_j \in E(G)\}$ and $N_{i}[v_i] = N_{i}(v_i) \cup \{v_i\}$. Further, we define, $N_{i}^{2}(v_i) = \{v_j \mid j >i ~\text{and} ~d_{G}(v_i,v_j) \leq 2)\}$ and $N_{i}^{2}[v_i] = N_{i}^{2}(v_i) \cup \{v_i\}$.

For other notations and graph theoretic terminology, we refer~\cite{tdbook}. In this paper, we consider only simple and connected graphs with at least $3$ vertices. Also, for a positive integer $n$, we use the standard notation, $[n]=\{1,2,\ldots, n\}$.

\section{Strongly Chordal Graphs}
\label{sec:1}
Strongly chordal graphs is an important subclass of chordal graphs introduced by several researchers in the literature \cite{sc1,sc2,sc3}. Strongly chordal graphs includes interval graphs, block graphs, directed path graphs, and trees as subclass. Many variations of domination are polynomial-time solvable on strongly chordal graphs, see \cite{sc4,sc1,sc8,sc5,sc6,sc7}. There are many equivalent definitions of strongly chordal graphs. We follow, the definition given in \cite{sc2}.

Let $G=(V,E)$ be a graph. A vertex $v \in V$ is called \emph{simple} if the vertices in  the closed neighbourhood of $v$ can be ordered, $N_{G}[v]=\{v_1,v_2, \ldots, v_r\}$ where $v_1=v$, such that $N_{G}[v_i] \subseteq N_{G}[v_j]$ for $1 \leq i \leq j \leq r$. A graph $G$ is \emph{strongly chordal} if every induced subgraph of $G$ has a simple vertex. An ordering $\alpha = (v_1,v_2, \ldots v_n)$ of vertices of $V$ is called \emph{strong elimination ordering}(SEO) if $v_j ,v_k \in N_{i}[v_i]$ implies $N_{i}[v_j] \subseteq N_{i}[v_k]$ for $1 \leq i \leq j \leq k \leq n$.

\begin{figure}[ht]
\begin{center}
\includegraphics[width=0.4\textwidth]{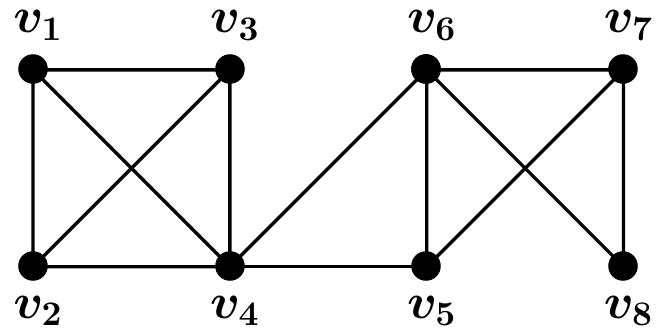}
\caption{An example of a strongly chordal graph.}
\label{fig:2}
\end{center}
\end{figure}

Many algorithms are studied to recognise a strongly chordal graph, $G=(V,E)$. In \cite{seo1,sc3}, the authors designed an $O(|V|^3)$-time algorithm to recognise a strongly chordal graph. In \cite{seo2}, a $O(|E| (\text{log}|E|)^{2})$-time algorithm is given which later, improved to $O(|E| \text{log}|E|)$-time algorithm in \cite{seo3}. Spinrad in \cite{seo4}, gave an $O(n^2)$-time algorithm to recognise a strongly chordal graph. The same algorithm also computes, a strong elimination ordering, if the graph is strongly chordal. The graph in Fig.~\ref{fig:2}, is a strongly chordal graph with strong elimination ordering $\alpha= (v_1,v_2,v_3,v_4,v_5,v_6,v_7,v_8)$.

Given a strongly chordal graph $G=(V,E)$ and a strong elimination ordering $\alpha = (v_1, v_2, \ldots v_n)$ of the vertex set $V$, we have the following observation.
\begin{obs}
\label{obs1}
If $v_iv_j \in E(G)$ such that $i <j$, then $N_{i}[v_i] \subseteq N_{i}[v_j]$. Further, if $F(v_i)=v_k$ then $N_{i}[v_j] \subseteq N_{i}[v_k]$.
\end{obs}
%\begin{proof}
%Proof directly follows using strong elimination ordering.
%\qed
%\end{proof}

Let $G=(V,E)$ be a strongly chordal graph and $\alpha = (v_1,v_2, \ldots v_n)$ be  its SEO. For a vertex $v_i \in V$, $F(v_i)$ denotes the highest index neighbour of $v_i$ according to SEO, where $i<n$. In particular $F(v_n)=v_n$. Our algorithm is an iterative algorithm which process the vertices as they appear in SEO. We use the following labels on the vertices during the execution of the algorithm to construct a minimum semi-TD-set of $G$.

\noindent
$
  D(v_i) =
  \begin{cases}
     0 & \text{if $v_i$ is not dominated}, \\
     1 & \text{if $v_i$ is dominated}. \\
  \end{cases}
$

\noindent
$
  L(v_i) =
  \begin{cases}
     0 & \text{if $v_i$ is not selected}, \\
     1 & \text{if $v_i$ is selected but no vertex $u$ is selected such that $d_{G}(u,v_i)\leq 2$}, \\
     2 & \text{if $v_i$ is selected and a vertex $u$ is also selected such that $d_{G}(u,v_i)\leq 2$}.\\
  \end{cases}
$

\noindent
$ m(v_i) =
  \begin{cases}
     k & {\text{if vertex $v_k \in N_{G}[v_i]$ is selected but a vertex $u$ need to be selected such that $d_{G}(u,v_k)\leq 2$}}, \\
     0 & \text{otherwise}. \\
  \end{cases}
$
\medskip

Further, $B_{i}[v_i]$ represents the set of neighbours $v_k$ of $v_i$ such that $N_{i}[F(v_i)] \subseteq N_{i}[v_k]$ and one of the neighbour of $v_k$ is already dominated. Formally, $B_{i}[v_i]= \{v_k \in N_{i}[v_i] \mid N_{i}[F(v_i)] \subseteq N_{i}[v_k]$ and there is a vertex $w \in N_{G}[v_k]$ such that $D(w)=1\}$. We note that, if $v_k \in B_{i}[v_i]$ then $N_{i}[F(v_i)] \subseteq N_{i}[v_k]$. Also by Observation~\ref{obs1}, $N_{i}[v_k] \subseteq N_{i}[F(v_i)]$. Hence, for a vertex $v_k \in B_{i}[v_i]$, we have $N_{i}[v_k] = N_{i}[F(v_i)]$.   In our algorithm, we also use two special types of operation on a particular vertex $v_i$, MARK$(v_j)$ and UNMARK$(v_j)$ which are defined as follows: if $F(v_i)= v_j$ then the operation MARK$(v_j)$ updates $m(v_k)=i$ for all $v_k \in N_{G}[v_j]$ such that $k>i$. While in operation UNMARK$(v_j)$ we update $m(v_k)=0$ for all $v_k \in N_{G}[v_j]$. Before designing the algorithm, we first prove the following results.

\begin{figure}[ht]
\includegraphics[width=0.7\textwidth]{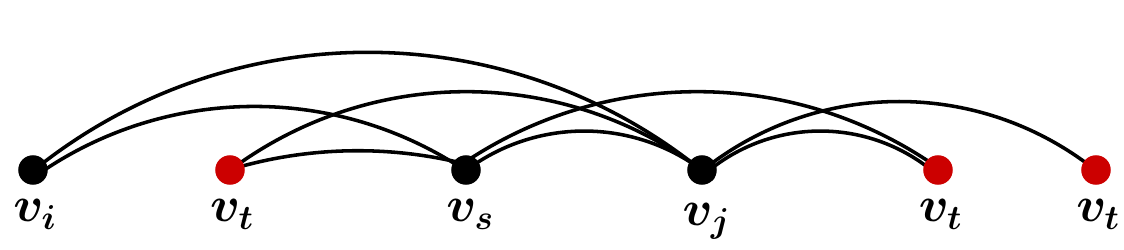}
\caption{Positions of $v_t$, when $d_{G}(v_i,v_t) = 2$ in Lemma~\ref{l:lemma1}}
\label{fig:3}
\end{figure}

\begin{lemma}
\label{l:lemma1}
For a vertex $v_i$, let $F(v_i)=v_j$ such that $i<j$ then $N_{i}^{2}[v_i] \subseteq N_{i}[v_j]$.
\end{lemma}
\pf
Since $v_iv_j \in E(G)$, $N_{i}[v_i] \subseteq N_{i}[v_j]$ using the property of SEO. Now, consider a vertex $v_t \in V(G)$ such that $d_{G}(v_i,v_t) =2$ and $t>i$. Let $P=v_iv_sv_t$ be a shortest path between $v_i$ and $v_t$. An illustration of possible positions of $v_t$ is given in Fig.~\ref{fig:3}. As $v_j$ is the highest index neighbour of $v_i$, we have $s \leq j$. Indeed, $s >i$, as if $s<i$ then using property of SEO we have $v_iv_t \in E(G)$, a contradiction. Hence, we have $i<s<j$, implying that $N_{i}[v_s] \subseteq N_{i}[v_j]$. Consequently, we have $v_tv_j \in E(G)$. Hence, the result follows.
\qed

\begin{lemma}
\label{l:lemma2}
If $v_iv_j \in E(G)$ such that $i<j$ then $N_{i}^{2}[v_i] \subseteq N_{i}^{2}[v_j]$.
\end{lemma}
\pf
Let $F(v_i)=v_k$ in the SEO. Clearly, $j \leq k$ and  using property of SEO, $ N_{i}[v_i] \subseteq N_{i}[v_j] \subseteq N_{i}[v_k]$. Let $v_t \in N_{i}^{2}[v_i]$ be an arbitrary vertex. If $v_iv_t \in E(G)$, then $v_jv_t \in E(G)$ as $N_{i}[v_i] \subseteq N_{i}[v_j]$. Now suppose, $d_{G}(v_i,v_t)=2$ and $P=v_iv_sv_t$ be a shortest path between $v_i$ and $v_t$. Clearly $s > i$. Since $N_{i}[v_i] \subseteq N_{i}[v_j]$, we have $v_sv_j \in E(G)$. Hence, $d_{G}(v_t,v_j) \leq 2$. Therefore, the lemma follows.
\qed

\begin{lemma}
\label{l:lemma3}
If $v_iv_j \in E(G)$ such that $i<j$ and $F(v_j)=v_k$ then $N_{i}^{2}[v_i] \subseteq N_{i}^{2}[v_k]$.
\end{lemma}
\pf
Using property of SEO, we note that $N_{i}[v_i] \subseteq N_{i}[v_j]$. Therefore, for any vertex $v_r \in N_{i}[v_i]$, we have $d_{G}(v_k,v_r) \leq 2$. Now consider a vertex $v_t$ such that $d_{G}(v_i,v_t) = 2$ and $t > i$. Let $P=v_iv_sv_t$ be a shortest path between $v_i$ and $v_t$. If $s<i$, then using the property of SEO, we have $v_iv_t \in E(G)$, a contradiction. Hence, $s>i$. Now, if $s \leq j$ then using property of SEO, $N_{i}[v_s] \subseteq N_{i}[v_j]$, implying that $d_{G}(v_k,v_t) \leq 2$. Further, if $s>j$ then $N_{i}[v_j] \subseteq N_{i}[v_s]$. This implies that $v_jv_s \in E(G)$ and hence, $v_sv_k \in E(G)$. Consequently, we have $d_{G}(v_k,v_t) \leq 2$. Therefore the lemma follows.
\qed

\begin{lemma}
\label{l:lemma4}
If for a vertex $v_i$, $F(v_i)=v_j \neq v_i$ and $v_k \in N_{i}[v_i]$ then $N_{i}^{2}[v_k] \subseteq N_{i}^{2}[v_j]$.
\end{lemma}
\pf
The proof directly follows from the property of strong elimination ordering. For completeness, suppose $v_r \in N_{i}^{2}[v_k]$. We note that $N_{i}[v_k] \subseteq N_{i}[v_j]$. Hence, if $v_kv_r \in E(G)$, Lemma follows. Now, assume that $d_{G}(v_k, v_r) =2$ where $r \geq i$. Let $P = v_rv_av_k$ be a shortest path of length two in $G$. As $k,r \geq i$, if $a<i$ then using property of SEO, we have $v_rv_k \in E(G)$, a contradiction. Hence, $a\geq i$. Further, using fact that, $N_{i}[v_k] \subseteq N_{i}[v_j]$, we have $v_av_j \in E(G)$. Consequently, we have  $d_{G}(v_r,v_j) \leq 2$.
\qed

\begin{figure}[ht]
\includegraphics[width=0.95\textwidth]{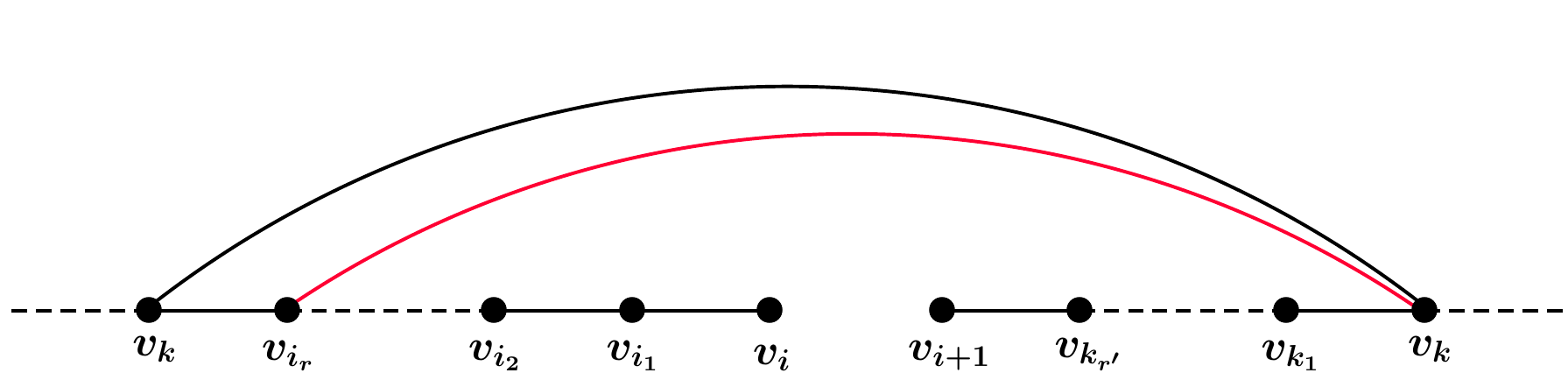}
\caption{An illustration of path in Lemma~\ref{l:lemma5}}
\label{fig:4}
\end{figure}

\begin{lemma}
\label{l:lemma5}
If for a vertex $v_i$, we have $F(v_i)=v_i$ then $i=n$.
\end{lemma}
\pf
On contrary, suppose $F(v_i)=v_i$ but $i \neq n$. Now since the graph is connected there exists a path joining the vertices $v_i$ and $v_{i+1}$. Suppose, $P$ is a shortest such path. Note that the path must contain a vertex $v_j$ such that $j<i$ and a vertex $v_k$ such that $k >i$. Let $v_j$ be the highest index vertex and $v_k$ be the least index in the path such that $v_jv_k \in E(G)$. Since, $P$ is a shortest path, we have $P = v_{i}v_{i_1}v_{i_2}\ldots v_{i_{r}}v_{j}v_{k}v_{k_1}v_{k_2}\ldots v_{k_{r'}}v_{i+1}$ where $j<i_r<i_{r-1}< \cdots < i_1 <i < i+1 < k_{r'} < \cdots v_{k_1} <v_{k}$. An illustration is given in Fig.~\ref{fig:4}. As $v_jv_k \in E(G)$, therefore, using the property of SEO, we have $v_{i_r}v_k \in E(G)$, a contradiction on choice of $P$. Hence, the result follows.
\qed

\begin{lemma}
\label{l:lemma6}
Let $G=(V,E)$ be a connected strongly chordal graph and $\alpha = (v_1,v_2, \ldots, v_n)$ be an SEO of $G$. For a vertex $v_i$, let $F(v_i)=v_j$. Suppose there is a path $P=v_{i}v_{i_1}v_{i_2} \ldots v_{i_{j}}v_{k}$ from $v_i$ to a vertex $v_k$ such that $k< i_1 < i_2 < \cdots < i_j < i $. If $v_k$ has a neighbour $v_s$ such that $s \geq i$ then $d_{G}(v_k,v_j) \leq 2$. Specifically, $v_iv_s,v_jv_s \in E(G)$.
\end{lemma}
\pf
The proof directly follows using property of strong elimination ordering. Hence, omitted.
\qed

\section{Algorithm for Semitotal Domination in Strongly Chordal Graphs}
\label{sec:2}

In this section, we propose, a linear-time algorithm to compute a minimum semi-TD-set in strongly chordal graphs. But, before designing the algorithm first, we discuss the idea of the algorithm.

\subsection*{\textit{\underline{Outline of the Algorithm}}}

Let $G=(V,E)$ be a strongly chordal graph and $\alpha = (v_1,v_2, \ldots ,v_n)$ be a strong elimination ordering of the vertex set of $G$. In our algorithm, we process the vertices iteratively as they appear in $\alpha$ and in each iteration we maintain a set $D_{st}$ containing the selected vertices. In the $i^{th}$-iteration we will process the vertex $v_i$. We ensure that the vertices having index at most $i-1$ are dominated by at least one vertex of $D_{st}$. Now in $i^{th}$-iteration, we update the set $D_{st}$ in the following way:
\begin{enumerate}
\item If the vertex $v_i$ is not dominated, and $v_i$ is not the last vertex, then we look for a vertex such that, $N_{i}[v_k] \subseteq N_{i}[F(v_j)]$ and there exists a vertex in $D_{st}$ which is at distance at most two from $v_k$.
\begin{itemize}
\item If such a vertex $v_k$ exists, we include it in $D_{st}$. We maintain the information that all the neighbours of $v_k$ are dominated now. Also, as there is already a vertex in $D_{st}$ which is at distance at most $2$ from $v_k$, we maintain the information that a distance two neighbour of $v_k$ is already selected. In the same iteration, we also check if there is a vertex $u \in D_{st}$, such that distance two neighbour of $u$ is not selected till $(i-1)^{th}$-iteration. If such a vertex exists, and $d_{G}(u,v_k) \leq 2$, then for all such vertices $u$, we update the information that a distance two neighbour of $u$ is selected, by updating $L(u)=2$ and UNMARK$(u)$-operation.

\item Otherwise, we include $F(v_i)$ in $D_{st}$ to dominate $v_i$ and updated the information that all the neighbours of $F(v_i)$ are dominated now. Further, using MARK$(F(v_i))$ operation, we also maintain the information that we need to select a distance two neighbour of $F(v_i)$, in one of the further iterations.
\end{itemize}

\item If the vertex $v_i$ is not dominated, and $v_i$ is the last vertex, that is, $v_i=v_n$, then we include $v_n$ in the $D_{st}$. We update the status of $v_n$ as dominated. Also, we may note that, all neighbours of $v_n$ are already dominated. Hence, there is a vertex in $D_{st}$ which is at distance at most $2$ from $v_n$. Therefore, we updated $L(v_n)=2$.

\item If $v_i$ is already dominated and $m(v_i)=0$, then we do not need to update $D_{st}$.

\item If $v_i$ is already dominated and $m(v_i) = k \neq 0$, that is, $v_{k}$ is selected in $D_{st}$ and a distance two neighbour of $v_k$ is still need to be selected then we updated $D_{st}$ as follows:
\begin{itemize}
\item If $i$ is less than the index of $F(v_k)$, we simply update $m(v_k)=0$.
\item If $v_i = F(v_k)$ and $F(v_i)=v_i$ then as $v_i$ is the last vertex, we include any neighbour of $v_i$ in $D_{st}$ and update $m(v_i)=0$.
\item If $v_i = F(v_k)$ and $F(v_i) \neq v_i$ then we include the highest index neighbour, $F(v_i)$ of $v_i$ in $D_{st}$. We update the information that all the neighbour of $F(v_i)$ are dominated now. We also update the information that a distance two neighbour of $v_k$ is selected by updating $L(v_k) =L(F(v_i)) =2$ and UNMARK$(v_k)$ operation. For all vertices, $u \in D_{st}$, such that distance two neighbour of $u$ is not selected till $(i-1)^{th}$-iteration and $d_{G}(u,v_k) \leq 2$, we update that a distance two neighbour of $u$ is selected, by updating $L(u)=2$ and UNMARK$(u)$-operation.
\end{itemize}
\end{enumerate}

\subsection*{\emph{\underline{Algorithm}}}

Now, we design a linear-time algorithm, Algorithm~\ref{l:algo}, to compute a minimum semi-TD-set in a strongly chordal graph, given its strong elimination ordering.

\begin{algorithm}[ht!]
\textbf{Input:} A Strongly Chordal graph $G=(V,E)$ and a SEO $ \alpha = (v_1,v_2, \ldots v_n)$ of $G$\\
\textbf{Output:} A minimum cardinality semi-TD-set $D_{st}$ of $G$.\\
\For{$i=1$ to $n$}{
	\If{$(D(v_i)=0$ and $F(v_i)\neq v_i)$}{
		Compute $B_{i}[v_i]$;\\
		\If {$(B_{i}[v_i]= \emptyset)$}{
			$L(F(v_i))=1$; // Let $F(v_i)=v_j$\\
			MARK$(v_j)$;\\
			$D(u)=1 $ $\forall$ $u \in N_{G}[v_j]$;
			}
		\Else{
			Let $k=$max$\{k' \mid v_{k'} \in B_{i}[v_i]\}$;\\
			$L(v_k)=2$;\\
			$D(u)=1 \forall u \in N[v_k]$;\\
			Let $C=\{v_r \mid v_r \in N[v_k]$ and $m(v_r) \neq 0\}$;\\
			\While{$(C \neq \emptyset)$}{
				Let $v_r \in C$ such that $m(v_r)=s$;\\
				$L(v_s)=2$;\\
				UNMARK$(v_s)$;\\
				$C=C\setminus\{v_r\}$;
				}
			}
		}
	\ElseIf{$(D(v_i) = 0$ and $F(v_i)=v_i)$}{
		$L(v_i)=2$;\\
		$D(v_i)=1$;
		}
	\ElseIf {$(D(v_i) = 1$ and $m(v_i) \neq 0)$}{
		Let $m(v_i)=s \neq 0$ and $F(v_s)=v_t$ where $t \geq s$;\\
		\If{$(i < t)$}{
			$m(v_i)=0$;
			}
		\ElseIf{$(i = t$ and $F(v_i) = v_i)$}{
			$L(u)=2$ for some $u \in N[v_i]$;\\
			$L(v_s)=2$ and $m(v_i)=0$;
			}
		\Else{
			$L(F(v_i))=2$;\\
			$L(v_s)=2$;\\
			$D(u)=1$ for all $u \in N[F(v_i)]$;\\
			Let $C=\{v_r \mid v_r \in N[F(v_i)]$ and $m(v_r) \neq 0\}$;\\
			\While{$(C \neq \emptyset)$}{
				Let $v_r \in C$ such that $m(v_r)=s'$;\\
				$L(v_s')=2$;\\
				UNMARK$(v_s')$;\\
				$C=C\setminus\{v_r\}$;
				}
			}
		}
}
Let $D_{st} = \{v_i \in V \mid L(v_i)= 2\} $;\\
return $D_{st}$;
\caption{\textbf{Minimum Semitotal Domination in Strongly Chordal Graphs}}
\label{l:algo}
\end{algorithm}

\newpage

Now to show the correctness of Algorithm~\ref{l:algo} and to prove that given a strongly chordal graph $G$ and its strong elimination ordering $\alpha$, it computes a minimum semi-TD-set of $G$, we prove the following lemmas (Lemma~\ref{l:lemma7} to Lemma~\ref{l:lemma12}).

\begin{lemma}
\label{l:lemma7}
At the beginning of the $i^{th}$-iteration, the following statements are true:
\begin{enumerate}
\item $D(v_j)=1$ for all $j \in [i-1]$.
\item $m(v_j)=0$ for all $j \in [i-1]$.
\item If $L(v_k) = 1$ for $k < i$, then there exists a neighbour $v_{k'}$ of $v_k$ such that $m(v_{k'})=k$ and $k' \geq i$.
\end{enumerate}
\end{lemma}
\pf
The proof directly follows from Algorithm~\ref{l:algo}.
\qed

\vspace{0.5cm}

For each $i \in [n]$, let $D_i = \{v_j \in V(G) \mid L(v_j) > 0\}$ be the set of selected vertices after the $i^{th}$-iteration of the algorithm. Indeed, $D_n = \{v_j \in V(G) \mid L(v_j) > 0\}$ is the set of vertices selected by the algorithm after all the vertices are processed. In particular, suppose $D_0 = \emptyset$. In order to prove the correctness of the algorithm first we will show that $D_n$ is a semi-TD-set of $G$ and then we will show that $D_i$ is contained in a minimum semi-TD-set of $G$ for each $i \in ([n] \cup \{0\})$.

\begin{lemma}
\label{l:lemma8}
$D_n$ is a semi-TD-set of strongly chordal graph $G$.
\end{lemma}
\pf
In the $i^{th}$-iteration of the algorithm, we process the vertex $v_i$. First we check whether the vertex $v_i$ is dominated or not using the vertices selected till $(i-1)^{th}$-iteration. Indeed, in the $i^{th}$-iteration, $D(v_i)=0$ represents that the vertex $v_i$ is not dominated by the current set $D_{i-1}$ and hence, we pick a neighbour of $v_i$ to dominate $v_i$. As we visit every vertex in some iteration, after the $n^{th}$-iteration, all vertices will be dominated by set $D_{n}$. Further in any iteration, if we pick any vertex $v_s$ such that no vertex $v_r$ is already selected such as $d_{G}(v_r,v_s) \leq 2$ then we perform MARK$(v_s)$ operation. Let $F(v_s)=v_t$. According to the algorithm, we select a vertex $u$ such $d_{G}(v_s,u) \leq 2$ in some $j^{th}$-iteration where $j \leq t \leq n$ and update $L(v_s)=2$. Hence, after the $n^{th}$-iteration we note that $L(v_i)=0 $ or $2$ for all $i \in [n]$. Hence $D_n$ is a semi-TD-set of $G$.
\qed

\medskip
Using Lemma~\ref{l:lemma8}, we note that the set $D_n$ is a semi-TD-set of the strongly chordal graph $G$. Now we claim that for each $i \in ([n] \cup \{0\})$, there is a minimum semi-TD-set $D'$ of $G$ containing $D_{i}$. We will prove this using induction on $i$. If $i=0$, $D_{0}= \emptyset$ is contained in any minimum semi-TD-set of $G$. Now we assume that there is minimum semi-TD-set $D$ of $G$ containing $D_{i-1}$. In the $i^{th}$ iteration, depending upon the several cases of the algorithm we prove the following lemmas to show that there is a minimum semi-TD-set $D'$ of $G$ containing $D_{i}$.

\begin{lemma}
\label{l:lemma9}
If $D(v_i)=0$, $F(v_i) \neq v_i$, and $B_{i}[v_i] =  \emptyset$ then there exists a minimum semi-TD-set $D'$ of $G$ such that $(D_{i-1} \cup \{F(v_i)\}) \subseteq D'$.
\end{lemma}
\pf
Using induction hypothesis, we have a minimum semi-TD-set $D$ of $G$ containing $D_{i-1}$. Clearly, if $F(v_i) \in D$ then $D'=D$ is the required minimum semi-TD-set. Hence, we assume that $F(v_i) \notin D$. Let $F(v_i)=v_j$. Since $F(v_i) \neq v_i$, $j>i$. Let $v_k$ and $v_l$ be the minimum index vertices in $D$ dominating $v_i$ and $v_j$ respectively. We note that $l,k \neq j$. Now we prove the result by considering the following two cases:

\vspace{0.9cm}

\noindent \textbf{Case 1:} If $v_k=v_l$

In this case $v_k$ is the least index vertex in $D$ dominating both $v_i$ and $v_j$. As $v_j= F(v_i)$ and $k \neq j$, therefore $k < j$. Indeed $v_k \notin D_{i-1}$ as if $v_k \in D_{i-1}$ then $D(v_i)=1$, a contradiction. Further, if $k<i$ then using Lemma~\ref{l:lemma7}, we have $D(v_k)=1$. Since $N_{i}[v_j] \subseteq N_{i}[v_j]$ and $v_kv_j \in E(G)$ with $D(v_k)=1$ hence, $v_j \in B_{i}[v_i]$, a contradiction as $B_{i}[v_i] = \emptyset$. Therefore $k \geq i$. As $v_k, v_j \in N_i[v_i]$ and $i \leq k <j$, using the property of SEO, we have $N_{i}[v_k] \subseteq N_{i}[v_j]$ and hence using Lemma~\ref{l:lemma7}, we note that $D' = (D \setminus \{v_k\}) \cup \{v_j\}$ is a dominating set of $G$ such that $|D'|=|D|$. Now in order to prove that $D'$ is a semi-TD-set, we need to prove the following claim.

\begin{claim}
\label{claim1}
For all $v_l \in N_{G}^{2}(v_k) \cap D$, there exists a vertex $v_r \in D'$ such that $d_{G}(v_l,v_r) \leq 2$ and there exists a vertex $u \in  N_{G}^{2}(v_j) \cap D'$.
\end{claim}
\pf
From previous arguments we note that $i \leq k <j$. Also $v_iv_k, v_kv_j, v_iv_j \in E(G)$. Let $v_l \in N_{G}^{2}(v_k) \cap D$ be an arbitrary vertex. If $(N_{G}^{2}(v_l) \cap D)\setminus \{v_k\} \neq \emptyset$, that is, there is a vertex $u$ other than $v_k$ in $D$ such that $d_{G}(v_l,u) \leq 2$, then the result follows. Hence, we assume that $N_{G}^{2}(v_l) \cap D = \{v_k\}$.
%Further, if $d_{G}(v_l,v_j) \leq 2$, the result follows therefore, we assume that $d_{G}(v_l,v_j) \geq 3$.

If $v_lv_k \in E(G)$ then we note that $d_{G}(v_l,v_j) \leq 2$ as $v_kv_j \in E(G)$, and the result follows. Assume that $d_{G}(v_k,v_l) = 2$. Let $v_lv_{l'}v_k$ be any shortest path between $v_l$ and $v_k$ in $G$. If $l' \geq i$ then using the property SEO, we have $v_{l'}v_j \in E(G)$, implying that $d_{G}(v_l,v_j) \leq 2$, hence the result follows.

Now suppose, $l'<i$. We have $v_{l'}v_l, v_{l'}v_{k} \in E(G)$. Therefore, if  $l' <l$, then using property of SEO, $v_lv_k \in E(G)$, contradicting the assumption that $d_{G}(v_k,v_l) = 2$. Hence, we have $l<l'<k$. Further, using Lemma~\ref{l:lemma4}, we note that if  $v_l$ has a neighbour $v_s$ such that $s \geq i$, then $d_{G}(v_l,v_j) \leq 2$. Hence, the result follows. Suppose, $v_l$ has no neighbour $v_s$ such that $s \geq i$. In this case, we show that either the vertex $v_l$ already has a distance two neighbour in $D$ or we can remove $v_l$ to get a semi-TD-set of smaller cardinality.

Here, first we show that there exists a vertex $v_{a'} \in D$ such that $a' \neq k$ and $d_{G}(v_j,v_{a'}) \leq 2$. We have $v_{l'}v_k \in E(G)$ and since $l'<i$, $D(v_{l'})=1$. We note that, $N_{i}[v_j] \nsubseteq N_{i}[v_k]$ as if $N_{i}[v_j] \subseteq N_{i}[v_k]$ then $v_k \in B_{i}[v_i]$, a contradiction. Hence, we have a vertex $v_a \in N_{i}[v_j]$ such that $v_av_k \notin E(G)$ that is, $v_a$ can not be dominated by $v_k$. Hence, there exists a vertex $v_{a'} \in D \setminus \{v_k\} \subseteq D'$ such that $v_av_{a'} \in E(G)$. Therefore, we have a vertex $v_{a'} \in D \setminus \{v_k\} \subseteq D'$ such that $d_{G}(v_{a'},v_j) \leq 2$.

Now, suppose $v_l \notin D_{i-1}$. Using Lemma~\ref{l:lemma7} and the fact that $v_l$ has no neighbour $v_s$ such that $s \geq i$, the set $D'' = D' \setminus \{v_{l}\}$ is a semi-TD-set of $G$ such that $|D''| < |D|$, a contradiction as $D$ is a minimum semi-TD-set of $G$. If $v_l \in D_{i-1}$ and $L(v_l)=1$ then by Lemma~\ref{l:lemma7}, there is a neighbour of $v_l$ having index greater than $i$ which is marked for $v_l$, a contradiction as $v_l$ has no neighbour having index greater than $i$. If $v_l \in D_{i-1}$ and $L(v_l)=2$, then there exists a vertex $v \in D_{i-1}$ such that $d_{G}(v_l,u) \leq 2$. Also we have a vertex $v_{a'} \in D'$ such that $d_{G}(v_j,v_{a'}) \leq 2$ and hence, the claim follows.
\qed

Consequently, in this case the result follows. Now we consider the other case.

\medskip
\noindent \textbf{Case 2:} If $v_k \neq v_l$

Since we have $D(v_i)=0$, $v_k \notin D_{i-1}$. If $k \geq i$, then using property of SEO, we note that $N_{i}[v_k] \subseteq N_{i}[v_j]$. Also if $k<i$, then $N_{k}[v_k] \subseteq N_{k}[v_i]$ implying that $N_{i}[v_k] \subseteq N_{i}[v_i] \subseteq N_{i}[v_j]$. Hence, the set $D'=D \setminus \{v_k\} \cup \{v_j\}$ is a dominating set of $G$. Also as $v_j$ is dominated by $v_l$ in $D$ and $v_l \in D \setminus \{v_k\} \subseteq D'$ hence, for $v_j \in D'$ we have $v_l \in D'$ such that $d_{G}(v_j,v_l) \leq 2$. Consequently, to prove $D'$ is a semi-TD-set of $G$ we need to prove the following claim.

\begin{claim}
For any vertex $v_r \in N_{G}^{2}(v_k)$, there exists a vertex $u$ such that $d_{G}(v_r,u) \leq 2$.
\end{claim}
\pf
Let $v_r \in N_{G}^{2}(v_k) \cap D$ be an arbitrary vertex. Note that if there exists a vertex $x \in (N_{G}^{2}(v_r) \cap D) \setminus \{v_k\}$ then the result follows. Hence we assume that $N_{G}^{2}(v_r) \cap D =\{v_k\}$. If $k \geq i$, then we can give similar arguments as we gave in Claim~\ref{claim1}, to show that there exists a vertex $v_s \in D'$ such that $d_{G}(v_r,v_s) \leq 2$ and hence, the result follows. Let $k < i$. First assume that $r > k$.
%Using the property of SEO, we note that $N_{k}[v_k] \subseteq N_{k}[v_i]$.
As $v_kv_i \in E(G)$, $k<i$, and $v_j=F(v_i)$, using Lemma~\ref{l:lemma3}, we have $N_{k}^{2}[v_k] \subseteq N_{k}^{2}[v_j]$. Therefore in this case, $d_{G}(v_r,v_j) \leq 2$ and the claim follows.

Now suppose, $r<k$. In this case, if $v_rv_k \in E(G)$ and $v_r$ has a neighbour $v_s$ such that $s>i$ then using Lemma~\ref{l:lemma6}, we have $d_{G}(v_r,v_j) \leq 2$. Now, assume that $d_{G}(v_r,v_k)=2$ and $P=v_rv_{r'}v_k$ be a shortest path joining $v_r$ and $v_k$. Note that $r'>r$, otherwise using the property of SEO, we have $v_rv_k \in E(G)$. Now if $v_r$ has a neighbour $v_s$ such that $s>i$ then using Lemma~\ref{l:lemma6}, again we have,
%property of SEO, $v_sv_{r'}, v_sv_k, v_iv_s \in E(G)$, hence $v_sv_j \in E(G)$ and $d_{G}(v_r,v_j) \leq 2$. So we may conclude that if $r<k$ and $v_r$ has a neighbour $v_s$ such that $s>i$, then
$d_{G}(v_r,v_j) \leq 2$ and we are done.

Suppose $v_r$ does not any neighbour $v_{s}$ such that $s>i$. Note that if $v_r \in D_{i-1}$ then $L(v_r) =2$. Indeed, if $L(v_r) =2$ then by Lemma~\ref{l:lemma7}, there exist a neighbour $v_s$ of $v_r$ such that $s \geq i$ and $m(v_s)=r$. Since, $v_r$ has no neighbour having index greater than $i$, we have a contradiction. Further, if $v_l \notin D_{i-1}$, then using Lemma~\ref{l:lemma7} and the fact that $v_l$ has no neighbour having index greater that $i$, the set $D'' = D' \setminus \{v_l\}$ is a semi-TD-set of $G$ such that $|D''| < |D|$, a contradiction, as $D$ is a minimum semi-TD-set of $G$. Hence, the claim follows.
\qed

As the lemma follows in both the cases, therefore, the result follows.
\qed

\begin{lemma}
\label{l:lemma10}
If $D(v_i)=0$, $F(v_i) \neq v_i$, and $v_k \in B_{i}[v_i] \neq  \emptyset$ where $k=$ max$\{k' \mid v_{k'} \in B_{i}[v_i]\}$ then there exists a minimum semi-TD-set $D'$ of $G$ such that $(D_{i-1} \cup \{v_k\}) \subseteq D'$.
\end{lemma}
\pf
If $v_k \in D$ then the result follows. Suppose $v_k \notin D$. Let $v_r \in D$ be the least index vertex dominating $v_i$. Since $D(v_i)=0$, $v_r \notin D_{i-1}$. Note that if $r<i$, then using the property of SEO, we have $N_{r}[v_r] \subseteq N_{r}[v_i]$ and hence, $N_i[v_r] \subseteq N_{i}[v_i] \subseteq N_{i}[v_j] \subseteq N_{i}[v_j]$. Also, if $r>i$ then using the property of SEO, we have $N_{i}[v_r] \subseteq N_{i}[v_j] \subseteq N_{i}[v_k]$. Hence, using Lemma~\ref{l:lemma7} and the fact that $N_{i}[v_r] \subseteq N_{i}[v_k]$, we note that the set $D'= (D \setminus \{v_r\}) \cup \{v_k\}$ is a dominating set of $G$. Further, we note that $v_k \in B_{i}[v_i]$. Hence, there exists a vertex $w \in N_{G}[v_k]$ such that $D(w)=1$. Consequently, there exists a vertex $u \in N_{G}^{2}(v_k) \cap D_{i-1}$ such that $d_{G}(v_k,u) \leq 2$.  Now to prove that $D'$ is a semi-TD-set of $G$ we need to prove the following claim:

\begin{claim}
For any vertex $v_a \in N_{G}^{2}(v_r) \cap D'$, there exists a vertex $u \in D'$ such that $d_{G}(v_a,u) \leq 2$.
\end{claim}
\pf
If there exists a vertex $x \in (N_{G}^{2}(v_a) \cap D) \setminus \{v_r\}$ then the result follows. Hence we assume that $N_{G}^{2}(v_a) \cap D =\{v_r\}$.

\noindent \textbf{Case 1:} If $r \geq i$.

Since, $v_iv_r \in E(G)$ and $r \geq i$, the property of SEO implies, $N_{i}[v_r] \subseteq N_{i}[v_j]$. Now, using the fact that $v_{k} \in B_{i}[v_i]$, we have $N_{i}[v_r] \subseteq N_{i}[v_j] \subseteq N_{i}[v_k]$.  Let $v_a \in N_{G}^{2}(v_r) \cap D'$ be an arbitrary vertex. Since, in this case, $v_rv_k \in E(G)$ therefore, if $v_av_r \in E(G)$, then we have  $d_{G}(v_a,v_k) \leq 2$ and result follows. So for further cases, we assume that $d_{G}(v_a,v_r)=2$ and $P=v_av_{a'}v_r$ is a shortest path between $v_a$ and $v_r$ in $G$.

Suppose $a \geq i$. Since $a,r \geq i$, we have $a'>i$. Indeed if $a'<i$ then using property of SEO, we observe that, $v_av_r \in E(G)$, a contradiction to our assumption that $d_{G}(v_a,v_r)=2$. Further, in this case, using the fact that $N_{i}[v_r] \subseteq N_{i}[v_k]$, we have $v_{a'}v_k \in E(G)$. This implies that $d_{G}(v_a,v_k) \leq 2$ and hence, the result follows.
%Since we have $N_{i}[v_r] \subseteq N_{i}[v_j] \subseteq N_{i}[v_k]$, we have $N_{i}^{2}[v_r] \subseteq N_{i}[v_j] \subseteq N_{i}[v_k]$. Hence, $d_{G}(v_a,v_k) \leq 2$ and we are done.

Now suppose, $a<i$.
%If $v_av_r \in E(G)$ then $d_{G}(v_a,v_k) \leq 2$ and hence, the result follows. Suppose, $d_{G}(v_a,v_r) = 2$ and $v_{a}v_{a'}v_r$ be a shortest path in $G$.
Note that if $a'>i$ then using the the fact that $N_{i}[v_r] \subseteq N_{i}[v_k]$, we have $v_{a'}v_k \in E(G)$. Hence we have, $d_{G}(v_a,v_k) \leq 2$ and the result follows.
Now, assume that $a'<i$. In that case, we have $a<a'$. Indeed, if $a'>a$ then using property of SEO, we have $v_av_r \in E(G)$, a contradiction. Now if $v_a$ has a neighbour $v_s$ such that $s \geq i$, then suing property of SEO, $v_{a'}v_s,v_rv_s \in E(G)$. Indeed we have $v_sv_k \in E(G)$, hence, $d_{G}(v_k,v_a) \leq 2$. Now suppose $v_a$ has no neighbour $v_s$ such that $s \geq i$. Here, if $v_a \in D_{i-1}$ then $L(v_a)=2$ as if $L(v_a)=1$, then by Lemma~\ref{l:lemma7}, there exists a vertex $v_s$ where $s \geq i$ such that $v_av_s \in E(G)$ and $m(v_s)=a$, a contradiction. If $v_a \notin D_{i-1}$ then as $v_a$ has no neighbour having index greater than $i$, we note that the set $D''=D' \setminus \{v_a\}$ is a semi-TD-set of $G$, contradicting the choice of $D$.  And hence the claim follows.
\smallskip

\noindent \textbf{Case 2:} If $r < i$.

First suppose $a \geq i$. If $v_av_r \in E(G)$, then using the property of SEO, we note that $d_{G}(v_k,v_a) \leq 2$. Suppose, $d_{G}(v_a,v_r)=2$, and $P=v_av_{a'}v_r$ is a shortest path between $v_a$ and $v_r$ in $G$. Note that, if $a'<r$, then we have $v_rv_a \in E(G)$, contradiction to our assumption that $d_{G}(v_k,v_a) \leq 2$. Hence, $r'>k$. Also we have $v_rv_i \in E(G)$ therefore using property of SEO we have, $N_{r}[v_r] \subseteq N_{r}[v_i]$. If $a'<i$, then we have $v_iv_{a'} \in E(G)$ and $N_{a'}[v_{a'}] \subseteq N_{a'}[v_i]$, implying that $d_{G}(v_a,v_k) \leq 2$. Otherwise, we have $v_iv_{a'} \in E(G)$ with $a'>i$, implying, $v_kv_{a'} \in E(G)$, and hence, $d_{G}(v_a,v_k) \leq 2$. Therefore, in this case, the result follows.

Now, suppose $a <i$. If $v_a$ does not have a neighbour having index greater than $i$ then using Lemma~\ref{l:lemma7}, either $L(v_a)=2$ or $v_a \notin D_{i-1}$. In the former case we note that the set $D''=D' \setminus \{v_a\}$ is a semi-TD-set of $G$, contradicting the choice of $D$. If $v_a$ has a neighbour having index greater than $i$, then as we discussed in previous cases, using  the property of SEO, we may observe that the $d_{G}(v_k,v_a) \leq 2$. And hence the claim follows.
\qed

This proves the lemma.
\qed

\begin{lemma}
\label{l:lemma11}
If $D(v_i)=0$ and $F(v_i) = v_i$ then there exists a minimum semi-TD-set $D'$ of $G$ such that $(D_{i-1} \cup \{v_i\}) \subseteq D'$.
\end{lemma}
\pf
As $F(v_i) = v_i$, it follows from Lemma~\ref{l:lemma5} that $i=n$. Since the graph is connected, there exists a neighbour $v_r$ of $v_n$ such that $r<n$. Also, using Lemma~\ref{l:lemma7}, we note that $D(v_j)=1$ and $m(v_j)=0$ for $j \in [n-1]$, implying that,  $D(v_r)=1$. Indeed, $m(v_n)=0$ as if $m(v_n)=k \neq 0$, then $v_{k} \in D_{i-1}$. Also we have $v_kv_n \in E(G)$ hence, $D(v_n) \neq 0$, a contradiction. By induction hypothesis, we note that there exists a minimum semi-TD-set $D$ such that $D_{n-1} \subseteq D$. Suppose, $v_l \in D$ is the least index vertex dominating $v_n$. Clearly, $v_l \notin D_{i-1}$. Using Lemma~\ref{l:lemma7}, we note that the set $D'=(D\setminus \{v_l\}) \cup \{v_n\}$ is a dominating set of $G$, such that $|D| = |D'|$. Now in order to show that $D'$ is a semi-TD-set, we need to show that, for every vertex $v_{k} \in N_{G}^{2}(v_l) \cap D$ there exists a vertex $u$ such that $d_{G}(v_k,u) \leq 2$.

If there exists a vertex $x \in (N_{G}^{2}(v_k) \cap D) \setminus \{v_l\}$, then the result follows. Hence we assume that $N_{G}^{2}(v_k) \cap D =\{v_l\}$. If $v_k \notin D_{n-1}$ then we may note that $D'' = D' \setminus \{v_k\}$ is a semi-TD-set of smaller cardinality, a contradiction. Further, if $v_l \in D_{n-1}$ and $L(v_l)=1$ then by Lemma~\ref{l:lemma7}, there exists a vertex $v_{k'}$ such that $m(v_{k'})=k$ and $k' \geq n$, a contradiction. If $v_l \in D_{n-1}$ and $L(v_l)=2$, the result follows.  This proves the lemma.
\qed

\begin{lemma}
\label{l:lemma12}
If $D(v_i) \neq 0$, $m(v_i) = s \neq 0$ and $F(v_s)=v_t$ then the following hold:
\begin{enumerate}
\item If $i <t$ then there is a vertex $v_r \in N[v_s]$ such that $m(v_r)=s$ and $r (\neq i) > i$.
\item If $i=t$ and $F(v_i) = v_i$ then there is a minimum semi-TD-set $D'$ of $G$ such that $(D_{i-1} \cup \{u\}) \subseteq D'$ where $u \in N_{G}[v_i]$.
\item If $i = t$ and $F(v_i) \neq v_i$ then there is a minimum semi-TD-set $D'$ of $G$ such that $(D_{i-1} \cup \{F(v_i)\}) \subseteq D'$.
\end{enumerate}
\end{lemma}

\pf
Let $D(v_i) \neq 0$, $m(v_i) = s \neq 0$ and $F(v_s)=v_t$ where $t \geq s$. First we will prove the following claims are true before the start of $i^{th}$-iteration.

\begin{claim}
There is a vertex $v_a$ such that $N_G[v_a] \cap D_{i-1} = \{v_s\}$, where $a<s$ and $F(v_a)=v_s$.
\end{claim}
\pf
Suppose $v_s$ is selected in $a^{th}$-iteration of the algorithm. We note that the algorithm selects the vertex $v_s$ in $a^{th}$-iteration only if $D(v_a)=0$.  Furthermore, the algorithm updates $L(v_s)=1$ and performs the MARK$(v_s)$ operation only if $B_{a}[v_a]=\emptyset$. Hence, before the start of $a^{th}$-iteration, we have $D(v_a)=0$, $F(v_a)=v_s$ and $B_{a}[v_a]=\emptyset$. After the $a^{th}$-iteration the algorithm has selected $v_s$ to dominate $v_a$ and has performed MARK$(v_s)$ operation. Hence,  $v_s \in (N_G[v_a] \cap D_{i-1})$.  Further, if $s=a$ then $F(v_a)=v_a =v_s$. In this case, in the $a^{th}$-iteration, the algorithm would have updated $L(v_s)=2$, and no neighbour of $v_s$ would have marked, a contradiction. Hence $a<s$.

Now suppose there is another vertex $v_c \in (N_G[v_a] \cap D_{i-1})$ such that $c \neq s$ and suppose $v_c$ is selected in the $b^{th}$ iteration.
%Note that $v_cv_a \in E(G)$ and $F(v_a)=v_s$ hence, using the property of SEO, we have either $c<a$ or $v_cv_s \in E(G)$ and $c<s$.
In the $a^{th}$ iteration, we have $F(v_a)=v_s$ and $B_{a}[v_a] = \emptyset$. If $b<a$, that is, if $v_c$ is selected before the start of $a^{th}$-iteration, then the algorithm would have updated $D(v_a)=1$ in the $b^{th}$-iteration. This is a contradiction to the fact that $D(v_a) = 0$ in $a^{th}$ iteration. Now suppose that $b>a$. In this case we have $c>a$ and $v_cv_s \in E(G)$. Note that in the $a^{th}$-iteration the algorithm would have selected the vertex $v_s$ and have performed the MARK$(v_s)$ operation. Also we note that, $v_cv_s \in E(G)$. Now, in the $b^{th}$-iteration, after selecting the vertex $v_c$, the algorithm finds a vertex $F(v_a) = v_s \in N_{G}[c]$ such that $m(v_s)=s$. Hence, algorithm would have updated $L(v_s)=2$ and have performed the operation UNMARK$(v_s)$, a contradiction. This proves the claim.
\qed

\begin{claim}
\label{claim5}
There is no vertex $v_r \in D_{i-1}$ such that $d_{G}(v_r,v_s) \leq 2$
\end{claim}
\pf
On contrary, we assume that there exists a vertex $v_r \in D_{i-1}$ such that $d_{G}(v_r,v_s) \leq 2$. Let $v_s$ and $v_r$ were selected in the $x^{th}$ and $y^{th}$ iteration of the algorithm. First suppose, $x>y$. In this case, the algorithm selected $v_r$ first and then $v_s$ is selected. Since, $d_{G}(v_r,v_s) \leq 2$ either $v_rv_s \in E(G)$ or $v_r$ and $v_s$ have a common neighbour $v_a$ in $G$. We note that after the execution of $y^{th}$ iteration, the algorithm would have updated $D(v_a)=1$. And since $v_sv_a \in E(G)$, during the execution of $x^{th}$-iteration, the algorithm would have updated $L(v_s) = 2$, a contradiction.

Now we consider the case when $x<y$. In this case, the algorithm first selects $v_s$ to dominate $v_x$ and then $v_r$ to dominate $v_y$. Also note that since, $d_{G}(v_r,v_s) \leq 2$ either $v_rv_s \in E(G)$ or $v_r$ and $v_s$ have a common neighbour $v_a$ in $G$. Suppose $r<s$. If $v_rv_s \in E(G)$, then while selecting the vertex $v_r$, the algorithm would have updated $L(v_s)=2$, a contradiction. Now suppose $d_{G}(v_r,v_s) = 2$ and let $v_a$ be the common neighbour of $v_r$ and $v_s$. Note that if $a<r$ then using the property of SEO, we have $v_rv_s \in E(G)$, a contradiction. Hence $a>r$. Further as $L(v_s)=1$, in the $y^{th}$-iteration we have, $D(v_a)=1$ and $m(v_a)=s$. Hence, after the $y^{th}$-iteration, the algorithm would have updated $L(v_s)=2$, and performed UNMARK$(v_s)$ operation, a contradiction. This proves the claim.
\qed

Now we continue the proof of the lemma. From above two claims, we note that there exists a vertex $v_a$ such that $N_G[v_a] \cap D_{i-1} = \{v_s\}$, where $a<s$ and $F(v_a)=v_s$. Also, there is no vertex $v_r \in D_{i-1}$ such that $d_{G}(v_r,v_s) \leq 2$. We note that $D(v_i) \neq 0$. Also by Lemma~\ref{l:lemma7}, $D(v_j)=1$ for all $j \in [i-1]$. Further, there is atleast one vertex, specifically $v_t$, such that $m(v_t)=s$. Hence, if $i<t$, then the result follows.

Further, we will prove that the Lemma follows in other two cases as well. We note that, using induction hypothesis, we have a minimum semi-TD-set $D$ of $G$ such that $D_{i-1} \subseteq D$. Furthermore, $L(v_s)=1$ implies $v_s \in D_{i-1} \subseteq D$. As $D$ is a minimum semi-TD-set of $G$, there exists another vertex $v_k \in D$ such that $d_{G}(v_s,v_k) \leq 2$. Using Claim~\ref{claim5}, we note that $v_k \notin D_{i-1}$.

Now, first suppose $i=t$ and $F(v_i)=v_i$. We claim that there exists a minimum semi-TD-set of $G$ containing $D_{i-1} \cup \{u\}$ where $u \in N_{G}[v_i]$. Note that $F(v_i)=v_i$ implies $i=n$. Therefore, using Lemma~\ref{l:lemma7} and the fact that $D(v_i)=1$, we note that the set $D'= (D \setminus \{v_k\}) \cup \{u\}$ is a minimum semi-TD-set of $G$ where $u \in N_{G}[v_i]$.

Finally, suppose $i=t$ and $F(v_i) \neq v_i$. Now we need to show that there exists a minimum semi-TD-set containing $D_{i-1} \cup \{F(v_i)\}$.
%By induction hypothesis, there exists a minimum semi-TD-set $D$ of $G$ such that $D_{i-1} \subseteq D$. Furthermore, $L(v_s)=1$ implies $v_s \in D_{i-1} \subseteq D$. As $D$ is a minimum semi-TD-set of $G$, there exists another vertex $v_k \in D$ such that $d_{G}(v_s,v_k) \leq 2$. By claim $2.2$, we note that $v_k \notin D_{i-1}$.
In this case, if $k<s$ then using Lemma~\ref{l:lemma6}, we have $N_{i}[v_k] \subseteq N_{i}[F(v_i)]$. Further, if $k>s$ then using Lemma~\ref{l:lemma1}, we have $v_k \in N_{s}[v_i]$, and hence, $N_{i}[v_k] \subseteq N_{i}[F(v_i)]$. Therefore,
%we have either $k<s$ or $v_k \in N_{s}[v_s] \subseteq N_{s}[v_i]$.
in all the above cases, we have $N_{i}[v_k] \subseteq N_{i}[F(v_i)]$. Hence, the set $D'=D \setminus \{v_k\} \cup \{F(v_i)\}$ is a dominating set of $G$. Let $F(v_i) = v_j$. In order to show that $D'$ is a semi-TD-set of $G$, we need to prove the following claim.

\begin{claim}
For every vertex $v_r \in N_{G}^{2}(v_k) \cap D'$, there exists a vertex $u \in D'$ such that $d_{G}(v_r,u) \leq 2$.
\end{claim}
\pf
Let $v_r \in N_{G}^{2}(v_k) \cap D'$ be an arbitrary vertex. Here $N_{G}^{2}(v_r) \cap D = \{v_k\}$. Indeed, if there is a vertex $u$ such that $u \in (N_{G}^{2}(v_r) \cap D) \setminus  \{v_k\}$ such that $d_{G}(v_r,u) \leq 2$ then the result follows. We will prove the claim by considering the following two cases:

\medskip
\noindent \textbf{Case 1:} If $k <s$.

Note that the result follows if any of the three condition holds for a vertex $v_r \in N_{G}^{2}(v_k) \cap D'$: $(i)$ $d_{G}(v_r,v_j) \leq 2$, $(ii)$ $d_{G}(v_r,v_s) \leq 2$, and $(iii)$ $v_r \in D_{i-1}$ and $L(v_r)=2$. Suppose none of the three conditions hold for the vertex $v_r$. We need to consider the following two sub-cases:
\medskip

\noindent \textbf{Case 1.1:} If $v_kv_s\in E(G)$.

In this case, using Lemma~\ref{l:lemma2}, we note that, $N_{k}^{2}[v_k] \subseteq N_{k}^{2}[v_s]$. Hence, if $r>k$ then we have $d_{G}(v_s,v_r) \leq 2$, a contradiction. Suppose $r<k$. If $v_rv_k \in E(G)$ then $d_{G}(v_s,v_r) \leq 2$, a contradiction. Next, assume that $d_{G}(v_k,v_r) = 2$. Let $v_rv_{k'}v_k$ be a shortest path between $v_r$ and $v_k$ in $G$. Note that $k'>r$, otherwise using the property of SEO, we have $v_rv_k \in E(G)$, a contradiction. Further, we note that if $v_r$ has a neighbour $v_{r'}$ such that $r' \geq s$, then using Lemma~\ref{l:lemma6}, we have, $d_{G}(v_s,v_r) \leq 2$. By our assumption, if $v_r \in D_{i-1}$ then $L(v_r)=1$, and by Lemma~\ref{l:lemma7}, there exists a neighbour $v_{r'}$ of $v_r$ such that $r' \geq i \geq s$, implying that $d_{G}(v_s,v_r) \leq 2$, a contradiction. Therefore, $v_r \notin D_{i-1}$. Hence, using Lemma~\ref{l:lemma7}, and the fact that $v_r \notin D_{i-1}$, the set $D''=D' \setminus \{v_r\}$ is a semi-TD-set of smaller cardinality, a contradiction.
\medskip

\noindent \textbf{Case 1.2:} If $d_{G}(v_k,v_s) = 2$.

Let $v_kv_{k'}v_s$ be a shortest path between $v_k$ and $v_s$ in $G$. Here, if $k'<k$, then using the property of SEO, we have $v_kv_s \in E(G)$, a contradiction to our assumption that $d_{G}(v_k,v_s) = 2$. Hence, $k'>k$.

\noindent \textbf{Case 1.2.1:} If $r>k$.

As we have $v_kv_{k'} \in E(G)$ and $k<k'$, hence using Lemma~\ref{l:lemma2}, $N_{k}^{2}[v_k] \subseteq N_{k}^{2}[v_{k'}]$. Consequently, $v_r \in N_{k}^{2}[v_k]$ implies $v_r \in N_{k}^{2}[v_{k'}]$. If $v_rv_{k'} \in E(G)$, then we have $d_{G}(v_s,v_r) \leq 2$, a contradiction. Now, suppose $d_{G}(v_{k'},v_r) = 2$ and $P=v_{k'}v_{r'}v_r$ be a shortest path between $v_{k'}$ and $v_r$.

If $k'>s$, using property of SEO, we have $v_{k'}v_i \in E(G)$. Since, $F(v_s)=v_i=v_t$, $k'<i$. Also if $r' \geq i$, then using property of SEO, we have $v_{r'}v_j \in E(G)$. This implies, $d_{G}(v_j,v_r) \leq 2$, a contradiction. Hence, $r'<i$. Now we consider the following two cases: $(i)$ $r>k'$ and $(ii)$ $r<k'$. In the first case, we note that $r' > k'$ and hence $v_{r'}v_i \in E(G)$. Since, we have $r>k'$, $v_{k'}v_i, v_{r'}v_i \in E(G)$, using the property of SEO we have, $v_rv_i \in E(G)$ and hence $d_{G}(v_r,v_j) \leq 2$, a contradiction.

Now, consider the second case, $r<k'$. In this case, we may note that, if $v_r$ has a neighbour $v_a$ such that $a \geq i$, then using Lemma~\ref{l:lemma6}, we have $d_{G}(v_j,v_r) \leq 2$, contradiction. Otherwise, using Lemma~\ref{l:lemma7} and the fact that $v_r$ has no neighbour $v_a$ such that $a \geq i$, the set $D''=D' \setminus \{v_r\}$ is a minimum semi-TD-set of smaller cardinality, a contradiction.

Now, consider the case if $k' \leq s$. Since, $v_{k'}v_{s} \in E(G)$, using Lemma~\ref{l:lemma2}, we have $N_{k'}^{2}[v_{k'}] \subseteq N_{k'}^{2}[v_{s}]$. Therefore, if $r>k'$, then we have $d_{G}(v_r,v_s) \leq 2$, a contradiction. Suppose $r<k'$. We note that if $v_r$ has a neighbour $v_a$ such that $a \geq i$, then using Lemma~\ref{l:lemma6}, we have $d_{G}(v_s,v_r) \leq 2$, contradiction. Otherwise, using Lemma~\ref{l:lemma7} and the fact that $v_r$ has no neighbour $v_a$ such that $a \geq i$, the set $D''=D' \setminus \{v_r\}$ is a minimum semi-TD-set of smaller cardinality, a contradiction.

%Suppose $r>k$. In this case, if $v_kv_r \in E(G)$ then using property of SEO, $v_{k'}v_r \in E(G)$. Hence, we have $d_{G}(v_s,v_r) \leq 2$, a contradiction. Now, assume that $d_{G}(v_k,v_r) = 2$. Let $v_kv_{r'}v_r$ be a shortest path between $v_k$ and $v_r$ in $G$.

%If $k' \geq s$ then using the property of SEO, $s$ $N_{k}[v_k] \subseteq N_{k}[v_k']$ and $N_{k}^{2}[v_k] \subseteq N_{k}[v_k']$. If there exists a vertex $v_{r}$ such that $d_{G}(v_k,v_r) \leq 2$ and $r>k$ then $v_{r} \in N_{k}[v_{k'}]$ and as $v_sv_{k'} \in E(G)$ and hence, $d_{G}(v_s,v_r) \leq 2$, a contradiction.

\noindent \textbf{Case 1.2.2:} If $r<k$.

If $v_kv_r \in E(G)$ and $v_r$ has a neighbour $v_{r'}$ such that $r'>i$, then using Lemma~\ref{l:lemma6}, we have, $d_{G}(v_j,v_r) \leq 2$, a contradiction. Let $d_{G}(v_k,v_r)=2$ and $P=v_kv_{r'}v_r$ be a shortest path joining $v_k$ and $v_r$. Note that, $r'>r$. Here, if $v_r$ has a neighbour $v_{a}$ such that $a>i$, then using Lemma~\ref{l:lemma6}, we have, $d_{G}(v_j,v_r) \leq 2$, a contradiction. Hence, in any case we may assume that, $v_r$ has no neighbour $v_{r'}$ such that $r'>i$. Consequently, using Lemma~\ref{l:lemma7} and the fact that $v_r$ has no neighbour $v_a$ such that $a \geq i$, the set $D''=D' \setminus \{v_r\}$ is a minimum semi-TD-set of smaller cardinality, a contradiction. Hence, the result follows.

\medskip
\noindent \textbf{Case 2:} If $k > s$.

Since, $v_sv_i \in E(G)$ such that $s < i$ and $F(v_s)=v_t=v_i$, using Lemma~\ref{l:lemma1}, we have $N_{s}^{2}[v_s] \subseteq N_{s}[v_i]$. Hence, $v_kv_i \in E(G)$.
\medskip

\noindent \textbf{Case 2.1:} If $k<i$.

In this case, first suppose $r>k$. If $v_kv_r \in E(G)$, then using the property of SEO, we have $v_rv_i \in E(G)$ and $d_{G}(v_k,v_j) \leq 2$, a contradiction. Now, assume that $d_{G}(v_k,v_r) = 2$ and $P=v_kv_av_r$ be a shortest path between $v_k$ and $v_r$ in $G$. Note that if $a<k$, then using the property of SEO, we have, $v_kv_r\in E(G)$, a contradiction. Hence, $a>k$. If $a>i$, then using the property of SEO, we have $v_av_i \in E(G)$ and hence, $v_av_j \in E(G)$. Therefore, $d_{G}(v_k,v_j) \leq 2$, a contradiction. Now if $a<i$, then using the property of SEO, we have $v_av_i \in E(G)$. Further, as $r>k$ and $v_av_r \in E(G)$, hence using the property of SEO, $v_iv_r \in E(G)$, implying that $\rightarrow$ $d_{G}(v_k,v_j) \leq 2$, a contradiction.

Next, suppose $r<k$. If $v_kv_r \in E(G)$ and $v_r$ has a neighbour $v_{r'}$ such that $r' \geq i$ then using Lemma~\ref{l:lemma6}, we have $d_{G}(v_r,v_j) \leq 2$, a contradiction. Now we consider the case if $d_{G}(v_k,v_{r}) =2$. Let $v_rv_{k'}v_k$ be any shortest path between $v_k$ and $v_k$ in $G$. Here also, we observe that if $v_r$ has a neighbour $v_{r'}$ such that $r' \geq i$ then using Lemma~\ref{l:lemma6}, we have $d_{G}(v_r,v_j) \leq 2$, a contradiction. Hence, in both cases, $v_r$ has no neighbour $v_{r'}$ such that $r' \geq i$. Consequently, using this fact and Lemma~\ref{l:lemma7} we note that the set $D''=D' \setminus \{v_r\}$ is a semi-TD-set of smaller cardinality, a contradiction. Hence, the result follows.
\medskip

\noindent \textbf{Case 2.2:} If $k \geq i$.

Using Lemma~\ref{l:lemma4}, we note that $N_{i}^{2}[v_k] \subseteq N_{i}^{2}[v_j]$. Hence, if $r>i$ then $d_{G}(v_j,v_r) \leq 2$, a contradiction. Now consider the case when $r<i$. Similar to the previous cases, we note that if $v_r$ has  has neighbour having index greater than $i$ then $d_{G}(v_r,v_j) \leq 2$, a contradiction. Further, if $v_r$ has no neighbour having index greater than $i$, then using Lemma~\ref{l:lemma7} we note that, the set $D''=D \setminus \{v_r\}$ is a semi-TD-set of smaller cardinality, a contradiction. Hence, the result follows.
\qed

This proves the lemma.
\qed

\vspace{0.5cm}
Using Lemma~\ref{l:lemma7} to Lemma~\ref{l:lemma12}, we may conclude that the set $D_{n}$ is a minimum semi-TD-set of $G$. Further, we may note that the Algorithm~\ref{l:algo} can be implemented in linear-time. Hence, we have the following result.

\begin{theorem}
Given a strongly chordal graph $G=(V,E)$, with strong elimination ordering $\alpha$ of vertex set $V$, a minimum semi-TD-set of $G$ can be computed in linear time.
\end{theorem}

\section{Conclusion}
\label{con}
In this paper, we resolved the complexity status of \textsc{Minimum Semitotal Domination} problem in strongly chordal graph, which is an important subclass of chordal graphs. The complexity status of the problem in dually chordal graphs and tolerance graphs is still unknown. It would be interesting to investigate the complexity of the problem in these graph classes. Further, as the problem is NP-complete in planar graphs, designing approximation algorithms for the problem in planar graphs is a good research direction.

%\bibliography{refer1}

%\end{thebibliography}
\end{document}